\documentclass[11pt]{amsart}
\usepackage[a4paper,margin=1in]{geometry}
\usepackage{amsmath,amssymb,amsthm, mathtools}
\usepackage[dvipsnames]{xcolor}
\usepackage{tikz} \usepackage{pgfplots}
\usepackage[hidelinks]{hyperref}
\usepackage{enumitem} 

\newtheorem{theorem}{Theorem}
\newtheorem{lemma}{Lemma}

\newtheorem{definition}[theorem]{Definition}
\newtheorem{introtheorem}{Theorem}

\numberwithin{equation}{section} 

\newcommand{\R}{\mathbb{R}}	

\newcommand{\norm}[1]{\left\lVert#1\right\rVert}
\newcommand{\abs}[1]{\lvert#1\rvert}

\newcommand{\bvar}{\beta} 
\newcommand{\bi}{B} 
\newcommand{\tempo}{T} 
\newcommand{\dd}{\,{\rm d}}                 
\DeclareMathOperator{\sgn}{sgn}       

\title{Carnot geodesics with H\"older controls}
\author{Nicola Cavallucci}
\author{Enrico Le Donne}
\author{Sebastiano Nicolussi Golo}
\date{21 October 2026}
\thanks{All authors were partially supported by the Swiss National Science Foundation (grant 200021\_204501, \emph{Regularity of sub-Riemannian geodesics and applications}).
The example of Section~\ref{sec:example} and a first version of its proof were obtained with ChatGPT Astra-Ultra (OpenAI); the exposition was revised with the assistance of GPT-6 (Codex) and Claude Opus 5.5 (Anthropic).
The authors checked all arguments and take full responsibility for the content.
We thank the University of Fribourg for IT support.}
\begin{document}
\begin{abstract}
{
We construct an explicit length-minimizing curve in a polynomial rank-two sub-Riemannian structure on $\R^4$.
Its arclength parametrization is $C^{1,1/2}$ and is not $C^{1,\alpha}$ near an interior point for any $\alpha>1/2$.
Minimality is proved by an approximate calibration chosen separately for each competitor.
The example admits a minimizing lift to a rank-two Carnot group of step five. Generalizations are announced. 
}
\end{abstract}
\maketitle

 \tableofcontents

\section{Introduction}\label{sec:introduction}

A sub-Riemannian length minimizer is a horizontal curve that minimizes length between its endpoints; all regularity assertions below refer to its arclength parametrization.
Determining the regularity of such curves is a central problem of the theory; see the survey~\cite{Monti}.
Smoothness is known to fail: Chitour, Jean, Monti, Rifford, Sacchelli, Sigalotti and Socionovo~\cite{CJMRSSS} constructed a $C^2$ but not $C^3$ minimizer with a boundary singularity, and Rossi, Schiavoni Piazza and Socionovo~\cite{RSS} obtained interior singularities and branching in real-analytic structures.
The general question of $C^1$ regularity remains open~\cite{RSS}.
Known restrictions include the exclusion of corners~\cite{Corners} and the existence of tangent lines~\cite{MPV}.
Hakavuori and Le Donne also showed that the projection of every blowup of a geodesic in a Carnot group of step at least two to the quotient by its last layer is a geodesic; consequently, in a sub-Riemannian manifold of step $s$, every $s$-times iterated tangent of a geodesic is a line~\cite[Theorem~1.1]{Blowups}.
Finally, length minimizers are $C^1$ for rank-two structures of step at most four~\cite{BCJPS}, and so are the minimizers of analytic three-dimensional structures~\cite{BDFPR}. 
{ Moreover, if the sub-Riemannian manifold is a Carnot group of rank 2 and step $\leq 4$, length minimizers are $C^\infty$~\cite{Leonardi_Monti}.}
Regarding the $L_p$-H\"older estimates of~\cite{LZreg}, these results concern integral regularity and do not by themselves give a pointwise H\"older modulus.
The example below has exact derivative exponent $1/2$ at an interior point, so it is $C^1$ but not $C^2$ in arclength.

In this article we give explicit minimizing geodesics whose derivative has a prescribed H\"older exponent at an interior point.
The basic example is four-dimensional.

\begin{introtheorem}\label{thm:intro-A}
On $\R^4$ with coordinates $(x,y,z,w)$, let $P\coloneq z-\frac25xy$, and consider the sub-Riemannian structure for which the vector fields $X=\partial_x$ and $Y=\partial_y+x\partial_z+P^2\partial_w$ are orthonormal. 
Then, the curve $\gamma(t)=\bigl(t\sqrt{\abs t},\,t,\,\frac25\abs t^{5/2},\,0\bigr)$, $\abs t\le10^{-8}$, minimizes the sub-Riemannian length between its endpoints.
Its controls $(\frac32\abs t^{1/2},1)$ are $\frac12$-H\"older and not differentiable at $t=0$.
The arclength parametrization of $\gamma$ is $C^{1,1/2}$, and it is not $C^{1,\alpha}$ near $t=0$ for any $\alpha>\frac12$; in particular, $\gamma$ is not $C^2$.
\end{introtheorem}

The minimality asserted in Theorem~\ref{thm:intro-A} is proved in Theorem~\ref{thm:half-holder} at page \pageref{thm:half-holder}.
The fields $X, Y$ generate an eight-dimensional stratified Lie algebra of rank two and step five.
The corresponding Carnot group projects onto $(\R^4, \{X, Y\})$, preserving controls and their norms, so $\gamma$ lifts to a minimizing geodesic of a Carnot group, with the same regularity. 
Consequently, the result \cite{Leonardi_Monti} on the smoothness of geodesics in rank 2 and step 4 is sharp in the step.
The construction used in Theorem~\ref{thm:intro-A} extends to smaller exponents:

\begin{introtheorem}\label{thm:intro-B}
For every integer $n\ge0$, there is a rank-two Carnot group of step at most $2n+5$ containing a minimizing geodesic whose arclength parametrization is $C^{1,\beta}$, with $\beta=1/(2(n+1))$, and is not $C^{1,\alpha}$ near an interior point for any $\alpha>\beta$.
\end{introtheorem}

Products of these examples give, for every finite set of times and every choice of exponents of this form, a minimizing geodesic in a Carnot group of rank $2k$ whose velocity fails to be differentiable exactly at the prescribed $k$ times, with the prescribed exponents.
{ The Carnot lift, Theorem~\ref{thm:intro-B}, and the product statement will be presented in a companion article.}

{
A non-$C^2$ minimizing geodesic also occurs in a polynomial metabelian structure on $\R^3$.
Declare the vector fields
 $\widetilde X\coloneq\partial_x$ and $
 \widetilde Y\coloneq\partial_y+(x^2-y^4)^2\partial_z
$
orthonormal.
Their generated Lie algebra is metabelian: all commutators are of the form $f(x,y)\partial_z$ and hence commute.
The curve
 $\widetilde\gamma(t)\coloneq(t|t|,t,0),$  for $ |t|\le 10^{-6},
$
minimizes length between its endpoints\footnote{
This example highlights a minor exception to~\cite[Corollary~1.3]{Soc}: while its $C^2$ assertion generally follows from the nonminimality analysis in~\cite[Theorem~1.2]{Soc}, the isolated case $(a,b)=(2,4)$ indeed admits a geodesic with a loss of regularity at the junction of two smooth branches.
}, and its arclength parametrization is $C^{1,1}$ but not $C^2$ at the interior point corresponding to $t=0$. 
For a competitor $\eta$ with the same endpoints and $\mu\coloneq\max|(x^2-y^4)\circ\eta|>0$, one can use the locally Lipschitz approximate calibration
\begin{equation*}\label{eq:intro-metabelian-calibration}
 \begin{aligned}
 \widetilde\Psi_\mu(x,y,z)
 &\coloneq F(x,y)-\frac{z}{8\mu^2}-2\int_0^y f_\mu(v)\dd v,\\
 F(x,y)&\coloneq\frac{2|y|x}{\sqrt{1+4y^2}}
       +\frac{y}{2\sqrt{1+4y^2}}
       +\frac14\operatorname{arsinh}(2y),
 \end{aligned}
\end{equation*}
where $f_\mu(v)\coloneq\min\{(32\mu^2)^{1/3},\,2\mu^2/v^4\}$ for $v\ne0$.
}

{
{\bf Method.}
Minimality is proved by an approximate calibration.
The underlying principle is the classical Hamilton--Jacobi verification argument:   if a $C^1$ function $F$ satisfies $\abs{\nabla_HF}\le1$ and $\nabla_HF$ equals the velocity of some horizontal curve parametrized by arc length, then the increment of $F$ between its endpoints is at most its length and the curve is length minimizing~\cite{HL,Fathi,ABB}. 
Our proof uses an approximate calibration chosen separately for each
competitor.
Fix a horizontal curve $\eta$ with the same endpoints as $\gamma$
and such that $L(\eta)\le L(\gamma)$.
Set $\mu\coloneq\max|P\circ\eta|$, which we may suppose positive. We construct a function $\Psi_\mu$ and denote by
$\Delta\Psi_\mu$ its increment between the common endpoints.
We prove the two estimates
\[
 L(\gamma)-5000\mu
 \le \Delta\Psi_\mu
 \le L(\eta)-\frac{10^8}{768}\mu.
\]
Since $10^8/768>5000$, these inequalities force
$L(\eta)>L(\gamma)$, contradicting our assumption.
The first estimate follows by evaluating $\Psi_\mu$ along $\gamma$:
the cutoff and correction terms change its endpoint increment
by at most $5000\mu$; see~\eqref{eq6ad897d8}.
For the second estimate, the corrections are chosen so that
\[
 |\nabla_H\Psi_\mu|
 \le 1-\frac{P^2}{64\mu^2}
 \qquad\text{along }\eta;
\]
see~\eqref{eq6ad91b1a}.
Thus the endpoint increment is smaller than $L(\eta)$ by at least
$\int_\eta P^2/(64\mu^2)$.
Because $P\circ\eta$ attains absolute value $\mu$ and has a controlled
Lipschitz constant, a peak estimate bounds this integral from below
by $(10^8/768)\mu$; see~\eqref{eq6ad897d7}.

}

{\bf Relation to previous work.}
The polynomial squared-constraint setting has antecedents in the work of Liu and Sussmann~\cite{LS}, and in the recent development~\cite{CJMRSSS}.
{  Retrospective framework references about calibrations are the historical reference \cite{HL} and Fathi's treatment \cite[Sections~4.1--4.3]{Fathi}. }
 Lifting procedures that produce Carnot examples are standard: they can be abstractly obtained by Palais's integrability theorem~\cite{Palais1957}, or directly, as in~\cite{RSS}.
We claim no priority for the calibration principle or for the elementary estimates above.
The features developed here are the competitor-dependent choice of vertical multiplier and cutoff, and the area coordinate $z$.
On curves where $P=0$, the relation $\dd z=x\,\dd y$ gives $y\,\dd x=\frac32x\,\dd y$, so a polynomial relation of low degree prescribes the fractional exponent of the two branches of $\gamma$.
For Theorem~\ref{thm:intro-B}, a hierarchy of signed primitives of $P$ improves the calibration error, and balancing its gain and loss exponents gives $\beta=1/(2(n+1))$.

Since $\gamma$ is $C^1$, its tangent at $t=0$ is unique and is a line, consistent with~\cite{MPV, Blowups}: the singularity of $\gamma$ is invisible to blowups, including iterated ones, and appears only in the H\"older modulus of its controls.

{\bf Organization.}
{
Section~\ref{sec:example} presents the four-dimensional example and establishes the properties of the structure and the candidate curve in~\S\ref{sub:subriemannian-observations}.
Section~\ref{sec:proof-four} proves minimality: after two elementary lemmas, we construct the cutoff functions in~\S\ref{subs6ad897c1}, build the calibrating family in~\S\ref{subs6ad897c9}, and conclude in~\S\ref{sub:review-minimality-four-dimensions}.
}

\section{A $C^{1,1/2}$ example}\label{sec:example}

On  $\R^4$, we use coordinates $(x,y,{z},{w})$, and consider the function
\[
P\coloneq P(x,y,{z},{w}) \coloneq  {z}-\frac25xy,\]
 and the vector fields
 \begin{equation}
 X\coloneq \partial_x,\qquad
 Y \coloneq  \partial_y+x\partial_{z}+P^2\partial_{w}.
 \label{eq:review-vector-fields}
\end{equation}

\begin{theorem}\label{thm:half-holder}
In the  sub-Riemannian structure on  $\R^{ 4}$ where $X$ and $Y$ are an orthonormal frame, the curve
\begin{equation}\label{curve_half_holder}
 \boxed{\quad
 t \in \left [-10^{-8}, 10^{-8} \right] \longmapsto
 \gamma(t)\coloneq
 \left(t\sqrt{|t|},\ t,\ \frac25|t|^{5/2},\ 0\right),
 \quad} 
\end{equation}
minimizes the length between its endpoints.
\end{theorem}
 See Figure~\ref{fig_gamma} for a visual representation.


\begin{figure}[htbp]
  \centering
\colorlet{HorizontalProjectionColor}{.}
\begin{tikzpicture}
\begin{axis}[
  width=6.5cm, height=6.5cm, axis equal image,
  axis lines=middle, xlabel={$x$}, ylabel={$y$},
  xmin=-1.12, xmax=1.12, ymin=-1.12, ymax=1.12,
  xtick={-1,0,1}, ytick={-1,0,1},
  axis line style={HorizontalProjectionColor},
  tick style={HorizontalProjectionColor},
  ticklabel style={text=HorizontalProjectionColor},
  label style={text=HorizontalProjectionColor},
  clip=false]
\addplot[HorizontalProjectionColor,very thick,
  domain=-1:1,samples=401,variable=\t]
  ({\t*sqrt(abs(\t))},{\t});
\end{axis}
\end{tikzpicture}
\caption{Horizontal projection $t\mapsto(t\sqrt{|t|},t)$, plotted for $-1\le t\le1$.}
\label{fig_gamma}
\end{figure}

\subsection{Properties of the structure and the candidate curve}\label{sub:subriemannian-observations}

\leavevmode{
In this subsection, we express horizontality and length in coordinates, then verify that $\gamma$ is horizontal and satisfies $P\circ\gamma=0$.
These computations identify its controls and give the length bound used to restrict possible shorter competitors.
}

Let $X,Y,P$ be as in~\eqref{eq:review-vector-fields}.
\subsubsection{The horizontal structure}

  We shall compute  the horizontal derivatives of the polynomial $P$:
\begin{equation}\label{der_P}
 XP=-\frac25y,\qquad YP=\frac35x.
\end{equation}
Indeed, these identities follow by differentiating $P={z}-\frac25xy$:
\begin{eqnarray*}
 XP
 &\stackrel{\eqref{eq:review-vector-fields}}{=}&
 \partial_x\left({z}-\frac25xy\right)\\
 &=&-\frac25y,\\[0.5ex]
 YP
 &\stackrel{\eqref{eq:review-vector-fields}}{=}&
 (\partial_y+x\partial_{z}+P^2\partial_{w})
 \left({z}-\frac25xy\right)\\
 &=&-\frac25x+x+P^2\cdot0\\
 &=&\frac35x.
\end{eqnarray*}

 \begin{lemma}\label{lemc_horizontality}
 In the sub-Riemannian structure given by \eqref{eq:review-vector-fields}, a curve $t\mapsto (x(t), y(t), {z}(t), {w}(t) )$ is horizontal if and only if
it satisfies
\begin{equation}\label{conditions horizontality}
 \dot {z}=x\dot y,\qquad \dot {w}=P^2\dot y.\end{equation}
Moreover, in that case, its speed is $\sqrt{\dot x^2+\dot y^2}$ and its length $L$ is
$$
 L=\int \sqrt{\dot x^2+\dot y^2}\dd t.$$
 \end{lemma}

%
%
%

\subsubsection{Some features of the candidate geodesic}

\begin{lemma}\label{lem6ad897bc_bis}
The curve $\gamma$ defined in \eqref{curve_half_holder} has the following properties:
\begin{enumerate}[label = (\roman*)]
\item  $\gamma$ is horizontal, 
  
\item  \begin{equation}\label{Pnullo}
P\circ\gamma\equiv 0
\end{equation}
  
\item   $L(\gamma)\le3T$, where $T \coloneq 10^{-8}$.
\end{enumerate}
\end{lemma}

\begin{proof}

The controls of $\gamma$ are 
\begin{equation}\label{controls gamma}
\dot x = \tfrac32|t|^{1/2}\quad  \text{ and }\quad \dot y=1.
\end{equation}
The function  $P$ vanishes along $\gamma$ because
$$ P(\gamma(t)) = {z}(t)-\frac25x(t)y(t)
=
\frac25|t|^{5/2}  -\frac25 t |t|^{1/2} t =0.
$$
Finally, we verify that
$$ \dot {z}=\frac52 \frac25|t|^{3/2}\sgn(t)  =t|t|^{1/2}= x= x\dot y,$$
and, since $P$ vanishes along $\gamma$,
$$ \dot {w}=0=P^2\dot y .$$
Thus \eqref{conditions horizontality} are satisfied.

Regarding the length $ L(\gamma)$, we bound:
\begin{eqnarray*}
 L(\gamma)&=&\int_{-T}^T\sqrt{\dot x^2+\dot y^2}\dd t\\
& \stackrel{\eqref{controls gamma}}{=}&
 \int_{-T}^T\sqrt{\tfrac94|t|+1}\dd t\\
  &\le& \int_{-T}^T\sqrt{\tfrac94 T+1 }\dd t\\
 &\le& 2 T   \sqrt{\tfrac94 T+1 }\\
 &\le&2 T \sqrt{\tfrac94} = 3T,
\end{eqnarray*}
where in the last inequality, we have used that $T$ is small enough that $\tfrac94 T+1\leq \tfrac94 $.
\end{proof}

%
%
%
%

\section{Proof of minimality in four dimensions}\label{sec:proof-four}
Before proving Theorem~\ref{thm:half-holder}, we establish two elementary geometric lemmas in Section~\ref{sub:two-geometric-lemmas}.
Afterward, in Section~\ref{subs6ad897c1}, we construct and analyze cutoff functions, which we assemble in Section~\ref{subs6ad897c9} to define the approximate calibrating family.

\subsection{Two geometric lemmas}\label{sub:two-geometric-lemmas}

\leavevmode{
In this subsection, we record two estimates that will be used in the calibration argument.
The first bounds the error along $\gamma$ caused by the cutoff; the second converts a decrease in one component of a unit covector into a quantitative norm deficit.
}

%
%
%
%
%
%
%
%
%

\subsubsection{Defect of a tilted unit covector}\label{subs6ad8974d}
{ The following lemma gives a bound on a defect: given a unit vector
\(
u_a:=\bigl(a,\sqrt{1-a^2}\bigr),
\)
and a (co)vector \(v_r:=(r,1)\),
the defect is, by definition, the quantity
\[
\norm{v_r}- u_a\cdot v_r
=\sqrt{1+r^2}-\bigl(ar+\sqrt{1-a^2}\bigr),
\]
which vanishes when \(u_a \) is allined with $v_r$. The assumption
\[
0\le a\le \frac{r}{\sqrt{1+r^2}}
\]
means that \(u_a\) is tilted from this aligned direction towards the vertical direction \((0,1)\).
In our future application, the vector \(v=(\frac32|t|^{1/2},1)\) is the control of the curve \(\gamma\), and \(u_a=(a_\mu(t),b_\mu(t))\) represents the horizontal gradient of a potential \(F_\mu\) along \(\gamma\). Thus the defect measures the loss in the calibration caused by the cutoff.
}

\begin{lemma}\label{lem6ad8974e}
Let $r\ge0$ and $a\in[0,\,r/\sqrt{1+r^2}]$.
Then
\begin{equation}\label{eq6ad8974f}
    0 \le \sqrt{1+r^2} - \Bigl(ar+\sqrt{1-a^2}\Bigr) \le \frac{r^2}{2} ,
\end{equation}
and the first inequality is an equality when $a=r/\sqrt{1+r^2}$.
\end{lemma}

\begin{proof}
Put $a_*\coloneq r/\sqrt{1+r^2}\in[0,1)$ and $\varphi(a)\coloneq ar+\sqrt{1-a^2}$ for $a\in[0,1)$.
Then $\varphi'(a)=r-a/\sqrt{1-a^2}$.
The function $a\mapsto a/\sqrt{1-a^2}$ is increasing on $[0,1)$ and takes the value $r$ at $a_*$, so $\varphi'\ge0$ on $[0,a_*]$.
Hence $\varphi$ is increasing on $[0,a_*]$, and for $a\in[0,a_*]$
\[
    1=\varphi(0)\le\varphi(a)\le\varphi(a_*)
    =\frac{r^2}{\sqrt{1+r^2}}+\frac1{\sqrt{1+r^2}}=\sqrt{1+r^2} .
\]
 The right inequality $\varphi(a)\le\sqrt{1+r^2}$ is the first inequality in~\eqref{eq6ad8974f}, with equality at $a=a_*$.
 The left inequality $\varphi(a)\ge1$ gives $\sqrt{1+r^2}-\varphi(a)\le\sqrt{1+r^2}-1\le r^2/2$, by $\sqrt{1+\tau}\le1+\tau/2$ for $\tau\ge0$.
\end{proof}


\subsubsection{Horizontal gradients slightly shorter than one}\label{subs6ad89750}

{
The following lemma estimates how much the norm of a unit vector decreases when we reduce its second component.
Starting from \(u_a=(a,b)\), with \(b\ge1/2\), we replace \(b\) by a nonnegative number \(d\le b-\sigma\), where \(\sigma\ge0\).
The lemma shows that the resulting vector has norm at most \(1-\sigma/4\).
In our application, these vectors are the reference covector  \((a_\mu,b_\mu) \) and \((a_\mu,\bvar_\mu)=\nabla_H\Psi_\mu\), evaluated along the competitor \(\eta\).
Thus the lemma converts the decrease in the second component into a quantitative norm deficit, which provides the gain in the calibration argument.}

\begin{lemma}\label{lem6ad89751}
Let $a,b,d,\sigma\in\R$ with $a^2+b^2=1$, $b\ge1/2$, $\sigma\ge0$ and $0\le d\le b-\sigma$.
Then
\begin{equation}\label{eq6ad89752}
    \sqrt{a^2+d^2}\le 1-\frac\sigma4 .
\end{equation}
\end{lemma}

\begin{proof}
Since $a^2=1-b^2$, we have $a^2+d^2=1-(b-d)(b+d)$.
Here $b-d\ge\sigma$ and $b+d\ge b\ge1/2$, so $a^2+d^2\le1-\sigma/2$.
Moreover $\sigma\le b-d\le b\le1$, so $1-\sigma/4\ge0$ and $(1-\sigma/4)^2=1-\sigma/2+\sigma^2/16\ge1-\sigma/2$.
Taking square roots gives~\eqref{eq6ad89752}.
\end{proof}

\subsection{Explicit $C^1$ cutoff functions}\label{subs6ad897c1}

\leavevmode{
We construct $C^1$ coefficients that remove the singularity of the calibration at $y=0$.
Their differential identity will cancel unwanted terms in the horizontal gradient, while their derivative bounds will control the correction used in the minimality proof.
}

Define $\chi:\R\to\R$ by
\begin{equation}\label{eq6ad897c2}
    \chi(r)\coloneq
    \begin{cases}
    0,&r\le1,\\
    3(r-1)^2-2(r-1)^3,&1<r<2,\\
    1,&r\ge2.
    \end{cases}
\end{equation}
On $(1,2)$, $\chi'(r)=6(r-1)(2-r)\in(0,3/2]$; at $r=1$ and $r=2$ the one-sided values and derivatives match.
So $\chi$ is $C^1$ and increasing, with $0\le\chi\le1$ and $0\le\chi'\le3/2$; see Figure~\ref{fig_chi} for a visual representation. 

\begin{figure}[htbp]
  \centering
\begin{tikzpicture}
\begin{axis}[width=12cm, height=5cm, xmin=0, xmax=3, ymin=-0.1, ymax=1.15,
  axis lines=left, xlabel={$r$}, ylabel={$\chi(r)$}, xtick={1,2}, ytick={0,1}, samples=300]
\addplot[very thick] coordinates {(0,0) (1,0)};
\addplot[very thick, domain=1:2] {3*(x-1)^2 - 2*(x-1)^3};
\addplot[very thick] coordinates {(2,1) (3,1)};
\node at (axis cs:0.5,0.08) {$0$};
\node at (axis cs:1.5,0.55) {$3(r-1)^2-2(r-1)^3$};
\node at (axis cs:2.5,1.08) {$1$};
\end{axis}
\end{tikzpicture}
\caption{The cutoff function $\chi(r)$, plotted for $0\le r\le3$.}
\label{fig_chi}
\end{figure}

Furthermore, for ${\mu}>0$ and $y\in\R$ define
\begin{equation}\label{eq6ad897c3}
    \check{a}(y)\coloneq\frac{\frac32\abs y^{1/2}}{\sqrt{1+\frac94\abs y}},\qquad
    a_{\mu}(y)\coloneq\chi(\abs y/\sqrt {\mu})\,\check{a}(y),\qquad
    {b}_{\mu}(y)\coloneq\sqrt{1-a_{\mu}(y)^2} ,
\end{equation}
See Figure~\ref{fig_cut} for visual representations. 
\begin{figure}[htbp]
  \centering
\begin{tikzpicture} \begin{axis}[ width=16cm, height=7cm, xmin=-1,xmax=1, ymin=-.3,ymax=1.1, samples=300, axis lines=middle, xlabel={$y$}, legend style={ at={(0.98,0.98)}, anchor=north east, draw=none }, domain=-1:1 ]
\addplot[blue,thick] {1.5*sqrt(abs(x))/sqrt(1+2.25*abs(x))}; \addlegendentry{$\check{a}(y)$}
\addplot[ red, thick ] { (abs(x)<=0.25) * 0 + (abs(x)>0.25 && abs(x)<0.5) * (3*(abs(x)/0.25-1)^2 -2*(abs(x)/0.25-1)^3) * (1.5*sqrt(abs(x))/sqrt(1+2.25*abs(x))) + (abs(x)>=0.5) * (1.5*sqrt(abs(x))/sqrt(1+2.25*abs(x))) }; \addlegendentry{$a_{\mu}(y)$}
\addplot[ green!60!black, thick ] { sqrt( 1- ( (abs(x)<=0.25) * 0 + (abs(x)>0.25 && abs(x)<0.5) * (3*(abs(x)/0.25-1)^2 -2*(abs(x)/0.25-1)^3) * (1.5*sqrt(abs(x))/sqrt(1+2.25*abs(x))) + (abs(x)>=0.5) * (1.5*sqrt(abs(x))/sqrt(1+2.25*abs(x))) )^2 ) }; \addlegendentry{${b}_{\mu}(y)$} \end{axis} \end{tikzpicture}
\caption{Plots for the functions $  \check{a}(y) ,    a_{\mu}(y) , $ and $   {b}_{\mu}(y)$
     with values ${\mu} = 1/16 $, plotted for $-1\le y\le 1$.}
\label{fig_cut}
\end{figure}
\\
Furthermore, let ${\bi}_{\mu}$ be the odd function with ${\bi}_{\mu}(0)\coloneq0$ and
\begin{equation}\label{eq6ad897c4}
    {\bi}_{\mu}(y)\coloneq-\frac52\,y^{-5/2}\int_0^yt^{3/2}\,a_{\mu}'(t)\dd t, \qquad \text{ for } y\in \R_+ .
\end{equation}

We next analyze the functions $\check{a},a_{\mu},{b}_{\mu}, {\bi}_{\mu}$ as defined above.
\begin{lemma}\label{lem6ad897c5}
For ${\mu}>0$, the functions $\check{a},a_{\mu},{b}_{\mu}, {\bi}_{\mu}$ defined in \eqref{eq6ad897c2} and \eqref{eq6ad897c4} satisfy the following properties:
\begin{enumerate}[label = (\roman*)]
\item $a_{\mu},{\bi}_{\mu},{b}_{\mu}\in C^1(\R)$;
\item $a_{\mu}={\bi}_{\mu}=0$ on $[-\sqrt {\mu},\sqrt {\mu}]$;
\item $a_{\mu}=\check{a}$ for $\abs y\ge2\sqrt {\mu}$;
\item $0\le a_{\mu}\le\check{a}<1$;
\item \begin{equation}\label{eq6ad897c68}
a_{\mu}^2+{b}_{\mu}^2=1;\end{equation}
\item
\begin{equation}\label{eq6ad897c69}
{b}_{\mu}(y)\ge \left(1+\frac94\abs y\right)^{-1/2}.\end{equation}
\item On $\R$,
\begin{equation}\label{eq6ad897c6}
    a_{\mu}'+{\bi}_{\mu}+\frac25y\,{\bi}_{\mu}'=0 .
\end{equation}
\item For $\abs y\ge\sqrt {\mu}$,
\begin{equation}\label{eq6ad897c7}
    \abs{a_{\mu}'(y)}\le6\abs y^{-1/2},\qquad\abs{{\bi}_{\mu}(y)}\le8\abs y^{-1/2},\qquad\abs{{\bi}_{\mu}'(y)}\le35\abs y^{-3/2} .
\end{equation}
\item
\begin{equation}\label{eq6ad897c8}
    \int_{\R}\abs{{\bi}_{\mu}'(y)}^2\dd y\le\frac{1225}{\mu} .
\end{equation}
\end{enumerate}
\end{lemma}

\begin{proof}
\emph{Items i--vi.}
By~\eqref{eq6ad897c2}, $\chi$ is $C^1$, takes values in $[0,1]$, vanishes on $(-\infty,1]$, and equals one on $[2,\infty)$.
By~\eqref{eq6ad897c3}, $\check{a}$ is even and smooth away from zero, while $a_\mu$ vanishes on $[-\sqrt\mu,\sqrt\mu]$.
Thus $a_\mu$ is even and $C^1$ on $\R$, and $a_\mu=\check{a}$ for $|y|\ge2\sqrt\mu$.
Also, $0\le a_\mu\le\check{a}<1$ by~\eqref{eq6ad897c3}.
Consequently, $b_\mu=\sqrt{1-a_\mu^2}$ is an even $C^1$ function and $a_\mu^2+b_\mu^2=1$.
Moreover,
\begin{eqnarray*}
 b_\mu(y)
 &\stackrel{\eqref{eq6ad897c3}, (iv)}{\ge}&\sqrt{1-\check{a}(y)^2}\\
 &\stackrel{\eqref{eq6ad897c3}}{=}&(1+\tfrac94|y|)^{-1/2}.
\end{eqnarray*}
Since $a_\mu'$ is continuous and vanishes on $[0,\sqrt\mu]$, the integral in~\eqref{eq6ad897c4} is $C^1$ for $y>0$ and vanishes for $0<y\le\sqrt\mu$.
Hence $\bi_\mu$ is $C^1$ on $(0,\infty)$ and zero on $(0,\sqrt\mu]$.
Its odd extension is therefore $C^1$ on $\R$ and vanishes on $[-\sqrt\mu,\sqrt\mu]$.
This proves Items i--vi.

\emph{Item vii.}
For $y>0$, multiplying~\eqref{eq6ad897c4} by $y^{5/2}$ and differentiating gives
\begin{eqnarray*}
 \frac52y^{3/2}\bi_\mu(y)+y^{5/2}\bi_\mu'(y)
 &\stackrel{\eqref{eq6ad897c4}}{=}&-\frac52y^{3/2}a_\mu'(y).
\end{eqnarray*}
Dividing by $\frac52y^{3/2}$ proves~\eqref{eq6ad897c6} for $y>0$.
Since $a_\mu'$ and $\bi_\mu$ are odd, and $\bi_\mu'$ is even, the left-hand side of~\eqref{eq6ad897c6} is odd, so the identity holds also for $y<0$.
All three terms vanish near zero, proving the identity on $\R$.

\emph{Item viii.}
For $y>0$, differentiation of~\eqref{eq6ad897c3} yields
\begin{eqnarray*}
 \check{a}'(y)
 &\stackrel{\eqref{eq6ad897c3}}{=}&
 \frac34y^{-1/2}(1+\tfrac94y)^{-3/2}
 \le\frac34y^{-1/2},\\
 a_\mu'(y)
 &\stackrel{\eqref{eq6ad897c3}}{=}&
 \frac1{\sqrt\mu}\chi'(y/\sqrt\mu)\check{a}(y)
 +\chi(y/\sqrt\mu)\check{a}'(y).
\end{eqnarray*}
The first term vanishes unless $\sqrt\mu<y<2\sqrt\mu$.
On this interval, $1/\sqrt\mu\le2/y$, while $0\le\chi'\le3/2$ by~\eqref{eq6ad897c2} and $\check{a}(y)\le\frac32\sqrt y$ by~\eqref{eq6ad897c3}.
It follows that, for every $y>0$,
\begin{eqnarray*}
 |a_\mu'(y)|
 &\stackrel{\text{preceding bounds}}{\le}&
 \frac2y\cdot\frac32\cdot\frac32\sqrt y
 +\frac34y^{-1/2}\\
 &=&\frac{21}4y^{-1/2}\le6y^{-1/2}.
\end{eqnarray*}
Using this bound inside~\eqref{eq6ad897c4}, we obtain, for $y>0$,
\begin{eqnarray*}
 |\bi_\mu(y)|
 &\stackrel{\eqref{eq6ad897c4}}{\le}&
 \frac52y^{-5/2}\int_0^yt^{3/2}|a_\mu'(t)|\dd t\\
 &\stackrel{\text{preceding bound}}{\le}&
 15y^{-5/2}\int_0^y t\dd t\\
 &=&\frac{15}2y^{-1/2}\le8y^{-1/2}.
\end{eqnarray*}
Therefore, for $y>0$,
\begin{eqnarray*}
 |\bi_\mu'(y)|
 &\stackrel{\eqref{eq6ad897c6}}{\le}&
 \frac5{2y}\bigl(|a_\mu'(y)|+|\bi_\mu(y)|\bigr)\\
 &\stackrel{\text{preceding bounds}}{\le}&
 \frac52\left(6+8\right)y^{-3/2}\\
 &=&
 35y^{-3/2}.
\end{eqnarray*}
The absolute values of $a_\mu'$, $\bi_\mu$, and $\bi_\mu'$ are even, so these estimates also hold for $y<0$.
In particular, they prove~\eqref{eq6ad897c7}.

\emph{Item ix.}
By Item ii and the continuity of $\bi_\mu'$, this derivative vanishes on $[-\sqrt\mu,\sqrt\mu]$.
Since $\bi_\mu'$ is even,
\begin{eqnarray*}
 \int_\R|\bi_\mu'(y)|^2\dd y
 &=&2\int_{\sqrt\mu}^\infty|\bi_\mu'(y)|^2\dd y\\
 &\stackrel{\eqref{eq6ad897c7}}{\le}&
 2\cdot35^2\int_{\sqrt\mu}^\infty y^{-3}\dd y\\
 &=&2\cdot35^2\cdot\frac1{2\mu}
 =\frac{1225}{\mu}.
\end{eqnarray*}
This proves~\eqref{eq6ad897c8}.
\end{proof}

\subsection{An explicit family of potentials}\label{subs6ad897c9}



{
We assemble the cutoff coefficients into a family of auxiliary potentials \((F_\mu)_{\mu>0}\) whose horizontal gradients along \(\gamma\) agree with its normalized control outside the cutoff region.
We compute these gradients and bound the endpoint deficit by \(9\mu/2\).
In the proof of the theorem, we will add correction terms to these potentials to obtain the approximate calibrations \(\Psi_\mu\).
These corrections provide a strict bound on the horizontal gradient along each competitor, yielding a gain that compensates for the losses in the endpoint increment.
}

Let $X,Y,P$ be as in~\eqref{eq:review-vector-fields} and $a_{\mu},{\bi}_{\mu},{b}_{\mu}$ as in~\ref{subs6ad897c1}.
{
With an abuse of notation, we shall see these functions as functions on $\R^4$ by precomposing them with the projection in the second variable, e.g., $a_{\mu}(x,y,z,w)\coloneq a_{\mu}(y)$.
}
\begin{definition}[Family of auxiliary potentials] 
For ${\mu}>0$, define the $C^1$ function $F_{\mu}:\R^4\to\R$,

\end{definition}
\begin{equation}\label{eq6ad897ca}
    F_{\mu}(x,y,{z},{w})\coloneq a_{\mu}(y)\,x+{\bi}_{\mu}(y)\,{z}+\int_0^y{b}_{\mu}(t)\dd t .
\end{equation}
For   $\gamma$ as in~\eqref{curve_half_holder}, the \emph{endpoint increment} of $F_{\mu}$ is denoted by $\Delta F_{\mu}\coloneq F_{\mu}(\gamma({\tempo}))-F_{\mu}(\gamma(-{\tempo}))$.

\begin{lemma}\label{lem6ad897cb}
Let ${\mu}>0$.
Then
\begin{equation}\label{eq6ad897cc}
    XF_{\mu}=a_{\mu},\qquad YF_{\mu}={b}_{\mu}+{\bi}_{\mu}'P ,
\end{equation}
and
\begin{equation}\label{eq6ad897cd}
    0\le L(\gamma)-\Delta F_{\mu}\le\frac92{\mu} .
\end{equation}
\end{lemma}

\begin{proof} 
The coefficient functions and their derivatives in the following computation are evaluated at $y$.
Since $F_\mu$ is independent of $w$ and $P=z-\frac25xy$, we have
\begin{eqnarray*}
 XF_\mu
 &\stackrel{\eqref{eq:review-vector-fields},\,\eqref{eq6ad897ca}}{=}&a_\mu,\\
 YF_\mu
 &\stackrel{\eqref{eq:review-vector-fields}}{=}&
 \partial_yF_\mu+x\partial_zF_\mu+P^2\partial_wF_\mu\\
 &\stackrel{\eqref{eq6ad897ca}}{=}&
 (a_\mu'+\bi_\mu)x+\bi_\mu'z+b_\mu\\
 &\stackrel{\eqref{eq6ad897c6}}{=}&
 b_\mu+\bi_\mu'\left(z-\frac25xy\right)\\
 &=&b_\mu+\bi_\mu'P.
\end{eqnarray*}
This proves~\eqref{eq6ad897cc}.

For $t\in\R$, put $r(t)\coloneq\frac32|t|^{1/2}$.
The chain rule along the horizontal curve $\gamma$ gives, for $t\in[-\tempo,\tempo]$,
\begin{eqnarray*}
 (F_\mu\circ\gamma)'(t)
 &\stackrel{\eqref{controls gamma}}{=}&
 r(t)(XF_\mu)(\gamma(t))+(YF_\mu)(\gamma(t))\\
 &\stackrel{\eqref{eq6ad897cc}}{=}&
 r(t)a_\mu(t)+b_\mu(t)+\bi_\mu'(t)P(\gamma(t))\\
 &\stackrel{\eqref{Pnullo}}{=}&
 r(t)a_\mu(t)+b_\mu(t).
\end{eqnarray*}
Here we also used $y(\gamma(t))=t$, which follows from~\eqref{curve_half_holder}.

Define the pointwise defect on $\R$ by
\[
 D_\mu(t)\coloneq\sqrt{1+r(t)^2}-\bigl(a_\mu(t)r(t)+b_\mu(t)\bigr).
\]
The cutoff definition~\eqref{eq6ad897c3} gives
\begin{eqnarray*}
 0\le a_\mu(t)
 &\stackrel{\eqref{eq6ad897c3}}{\le}&
 \check a(t)=\frac{r(t)}{\sqrt{1+r(t)^2}},\\
 b_\mu(t)
 &\stackrel{\eqref{eq6ad897c3}}{=}&\sqrt{1-a_\mu(t)^2}.
\end{eqnarray*}
Thus Lemma~\ref{lem6ad8974e}, applied with $r=r(t)$ and $a=a_\mu(t)$, yields
\begin{equation}\label{eq:calibration-defect-bound}
 0\stackrel{\eqref{eq6ad8974f}}{\le}D_\mu(t)
 \stackrel{\eqref{eq6ad8974f}}{\le}\frac{r(t)^2}{2}
 =\frac98|t|.
\end{equation}
If $|t|\ge2\sqrt\mu$, then~\eqref{eq6ad897c2}--\eqref{eq6ad897c3} give $a_\mu(t)=\check a(t)$.
The equality case of Lemma~\ref{lem6ad8974e} therefore gives
\begin{equation}\label{eq:calibration-defect-support}
 D_\mu(t)=0\qquad\text{if }|t|\ge2\sqrt\mu.
\end{equation}

By~\eqref{controls gamma}, the speed of $\gamma$ is $\sqrt{1+r(t)^2}$.
The fundamental theorem of calculus and the formula for $(F_\mu\circ\gamma)'$ therefore give
\[
 L(\gamma)-\Delta F_\mu=\int_{-\tempo}^{\tempo}D_\mu(t)\dd t\ge0.
\]
For the upper bound, set $E_\mu\coloneq[-\tempo,\tempo]\cap[-2\sqrt\mu,2\sqrt\mu]$.
Then
\begin{eqnarray*}
 L(\gamma)-\Delta F_\mu
 &=&\int_{-\tempo}^{\tempo}D_\mu(t)\dd t\\
 &\stackrel{\eqref{eq:calibration-defect-support}}{=}&
 \int_{E_\mu}D_\mu(t)\dd t\\
 &\stackrel{\eqref{eq:calibration-defect-bound}}{\le}&\frac98\int_{E_\mu}|t|\dd t\\
 &\stackrel{E_\mu\subseteq[-2\sqrt\mu,2\sqrt\mu]}{\le}&
 \frac98\int_{-2\sqrt\mu}^{2\sqrt\mu}|t|\dd t\\
 &=&\frac94\int_0^{2\sqrt\mu}t\dd t\\
 &=&\frac94\cdot\frac{(2\sqrt\mu)^2}{2}
 =\frac92\mu.
\end{eqnarray*}
This proves~\eqref{eq6ad897cd} for every $\mu>0$, including when $2\sqrt\mu>\tempo$.
\end{proof}

\subsection{Proof of minimality}\label{sub:review-minimality-four-dimensions}

\leavevmode{
We now prove Theorem~\ref{thm:half-holder}.
For each competitor, we will consider an almost calibration $\Psi_\mu$ constructed from the map $F_\mu$ defined in~\eqref{eq6ad897ca}.
}

\begin{proof}[Proof of Theorem~\ref{thm:half-holder}]
We continue using the notation \begin{equation}\label{def_tempo} \tempo\coloneq10^{-8}.\end{equation}

Let   $\gamma$ be as in~\eqref{curve_half_holder}.
Let $\eta$ be a competitor for $\gamma$, i.e.,
let $\eta:[0, L]\to\R^4$ be a horizontal curve  with the same endpoints as $\gamma$ and with $L(\eta) \le L(\gamma)$.

We may assume that $\eta: t\mapsto (x(t), y(t), {z}(t), {w}(t) )$  is parametrized by arclength. Thus, we have

  \begin{equation}\label{endpoint1}\eta(0)=\gamma(-{\tempo})  , \end{equation}

    \begin{equation}\label{endpoint2}\eta(L)=\gamma({\tempo})  ,\end{equation}

      \begin{equation}\label{bound_pbal}
    |  \dot x|, \, |\dot y | \leq 1\quad \text{ and } \quad \dot x^2+ \dot y^2=1 , \qquad \text{almost everywhere},
       \end{equation}

     \begin{equation}\label{bound_length}L \stackrel{\text{def}}{=}L(\eta) \le L(\gamma)
     \le
     3{\tempo} , \end{equation}
 where the last inequality comes from Lemma~\ref{lem6ad897bc_bis}.iii.

  We claim the following bound on the coordinates of $\eta$:
        \begin{equation}\label{apriori_bound_coord}
        \abs{x(t)},\abs{y(t)}\le4{\tempo} \qquad  \forall t\in[0,L] .
          \end{equation}
Indeed, for every $t\in[0,L]$, absolute continuity gives
\begin{eqnarray*}
 |x(t)|
 &=& \left|x(0)+\int_0^t\dot x(s)\dd s\right|\\
 &\le& |x(0)|+\int_0^t|\dot x(s)|\dd s\\
 &\stackrel{\eqref{endpoint1},\,\eqref{curve_half_holder}}{=}&
 {\tempo}^{3/2}+\int_0^t|\dot x(s)|\dd s\\
 &\stackrel{\eqref{bound_pbal}}{\le}& {\tempo}^{3/2}+t\\
 &\stackrel{0\le t\le L,\ {\tempo}\le1}{\le}& {\tempo}+L\\
 &\stackrel{\eqref{bound_length}}{\le}&4{\tempo}.
\end{eqnarray*}
Similarly,
\begin{eqnarray*}
 |y(t)|
 &\le& |y(0)|+\int_0^t|\dot y(s)|\dd s\\
 &\stackrel{\eqref{endpoint1},\,\eqref{curve_half_holder}}{=}&
 {\tempo}+\int_0^t|\dot y(s)|\dd s\\
 &\stackrel{\eqref{bound_pbal}}{\le}&{\tempo}+t\\
 &\stackrel{t\le L}{\le}&{\tempo}+L\\
 &\stackrel{\eqref{bound_length}}{\le}&4{\tempo}.
\end{eqnarray*}

We consider the real-valued curve
 \begin{equation}\label{def_p}
 p\coloneq P\circ\eta, \quad \text{ on } [0,L]
 \end{equation}
 and the constant (which will represent the cutoff parameter for the approximate calibrations)
 $${\mu} \coloneq \max\{|p(t)| \,:\, t\in [0,L]\}.$$

We claim that
$p$ is $3{\tempo}$-Lipschitz, i.e.,
        \begin{equation}\label{p_Lip}
        |p'| \leq 3{\tempo}.
          \end{equation}
Indeed, the composition $p=P\circ\eta$ is absolutely continuous, and the chain rule\footnote{We may identify the horizontal gradient $\nabla_H P=(XP,YP)$ with the horizontal vector field $ (XP)X+(YP)Y$, with the norm given by the sub-Riemannian structure.}, together with the preceding identities, almost everywhere give
\begin{eqnarray*}
 |p'|
 &=&
\abs{ \nabla_H P\cdot \dot\eta}\\
  &\stackrel{\text{Cauchy--Schwarz}}{\le}&  \| \nabla_H P\| \cdot \|\dot\eta\| \\
  &\stackrel{\eqref{bound_pbal}}{\le}&  \sqrt{ (X P)^2 + (Y P)^2 }   \cdot 1\\
 &\stackrel{\eqref{der_P}}{=}&
 \frac15\sqrt{4y(s)^2+9x(s)^2}\\
 &\stackrel{\eqref{apriori_bound_coord}}{\le}&
 \frac{4\sqrt{13}}5{\tempo}\\
 &
 {<}
 &3{\tempo}.
\end{eqnarray*}

We next claim that
        \begin{equation}\label{bound_M}
  {\mu} \leq 10^{-15}.
          \end{equation}
Indeed, for every $t\in[0,L]$,
\begin{eqnarray*}
 |p(t)|
 &=&\left| p(0)+ \int_0^t p'(s)\dd s\right|\\
  &\stackrel{\eqref{def_p}}{\leq} &\left| P(\eta(0))\right|+\left| \int_0^t p'(s)\dd s\right|\\
 &\stackrel{ \eqref{endpoint1}, \eqref{p_Lip}}{\le}&\left| P(\gamma(-{\tempo}))\right|+ 3{\tempo}t\\
 &\stackrel{\eqref{Pnullo}, t\le L}{\le}&0+ 3{\tempo}L\\
 &\stackrel{\eqref{bound_length}}{\le}&9{\tempo}^2\\
 &
 \stackrel{\eqref{def_tempo}}{ =}
&9\cdot10^{-16}<10^{-15}.
\end{eqnarray*}
Taking the maximum over $t\in[0,L]$ proves~\eqref{bound_M}.

Next, we perform a peak estimate on $|p|$, which is $3{\tempo}$-Lipschitz, vanishes at the boundary points, and has maximum equal to ${\mu}$. Namely,  the interval of
length ${\mu}/(3{\tempo})$ centered at any maximum point of $|p|$ must lie in $[0,L]$
and in this interval we have  $|p|\geq  {\mu}/2$. Consequently,
\begin{equation}\label{peak_estimate_used}
 \int_0^L p^2\dd s\ge\frac{{\mu}^3}{12 {\tempo} }.
\end{equation}

We distinguish the two cases ${\mu}>0$ and ${\mu}=0$.

\subsubsection{Case: ${\mu}>0$}
We suppose first that ${\mu}>0$.
In this case, we use the family $a_{\mu},{\bi}_{\mu},{b}_{\mu},F_{\mu}$ defined above associated to $\mu$ and we further consider the function
\begin{eqnarray}
\label{bvar1}
 {\bvar_{\mu}}&\coloneq& {b}_{\mu}+{\bi}_{\mu}'P-\frac{P^2}{8{\mu}^2}-4{\mu}^2\abs{{\bi}_{\mu}'}^2\\
\label{bvar2}
  &=&{b}_{\mu}-\frac{P^2}{16{\mu}^2}
    -\frac12\left(\frac{P}{2\sqrt2 {\mu}}-2\sqrt2 {\mu} {\bi}_{\mu}'\right)^2,
    \qquad\text{on }\R^4.
\end{eqnarray}

We consider the function $\Psi_{\mu}$ defined as
 \begin{equation}\label{def_Phi}
 \Psi_{\mu}\coloneq F_{\mu}-\frac{{w}}{8{\mu}^2}
          -4{\mu}^2\int_0^y|{\bi}_{\mu}'(t)|^2\dd t ,
\end{equation}
which is the (unique-up-to-constant)  $C^1$ function  with horizontal gradient satisfying
 \begin{equation}\label{def_Phi_gradient}
 \nabla_H\Psi_{\mu}= (a_{\mu},{\bvar_{\mu}}),\end{equation}.

 For  all $s\in[0,L]$ we consider the values 
 $$a\coloneq a_{\mu}(\eta(s)) ,\quad b\coloneq b_{\mu}(\eta(s)), \quad\bvar \coloneq {\bvar_{\mu}}(\eta(s)),\quad \text { and } \quad \sigma\coloneq \frac{p(s)^2}{2^4 {\mu}^2}.$$
We claim that these values satisfy
\begin{equation}\label{check_geometric_lemma}
\left\{\begin{array}{l}
{\bvar}\leq {b}-\sigma\\
\sigma\geq 0\\
a ^2+{b}^2=1\\
{b}\geq 1/2\\
{\bvar }\geq 0 .
\end{array}
\right.
\end{equation}
Indeed, 
the first comes from \eqref{bvar2}; the second one comes directly from the definition; the third one is \eqref{eq6ad897c68}.
Regarding the 4th equation,  by Lemma~\ref{lem6ad897c5}.6 and \eqref{apriori_bound_coord},
\begin{equation}\label{bmaggioremezzo}
 b \stackrel{\eqref{eq6ad897c69}}{\ge}(1+\tfrac94|y|)^{-1/2}
\stackrel{\eqref{apriori_bound_coord}}{ \ge}(1+9\tempo)^{-1/2}\stackrel{\eqref{def_tempo}}{ \ge}\tfrac12.
\end{equation}
Observe that by Lemma~\ref{lem6ad897c5}.2 we have $\bi_\mu'(y)=0$ for all $|y|\le\sqrt\mu$. Thus, 
\eqref{eq6ad897c7} gives the global bound
$|\bi_\mu'|\le35\mu^{-3/4}$.
Consequently, since $|P|\le\mu$ along $\eta$, the definition of $\bvar_\mu$ yields
\begin{eqnarray*}
 \bvar
 &\stackrel{\eqref{bvar1}}{\ge}& b-|\bi_\mu'|\,|P|
       -\frac{P^2}{8\mu^2}-4\mu^2|\bi_\mu'|^2\\
 & \stackrel{\eqref{bmaggioremezzo}}{\ge}& \frac12-35\mu^{1/4}-\frac18-4900\mu^{1/2}\\
 &\stackrel{\eqref{bound_M}}{\ge}&
 \frac38-35\cdot10^{-3}-4900\cdot10^{-6}>0.
\end{eqnarray*}

The inequalities \eqref{check_geometric_lemma} tell us that we can apply the geometric Lemma~\ref{lem6ad89751}, and obtain for  all $s\in[0,L]$ the bound
\begin{equation}\label{eq6ad91b1a}
    \abs{\nabla_H\Psi_{\mu}}(\eta(s))
    \stackrel{\eqref{def_Phi_gradient}}{=}
    \sqrt{a^2+{\bvar}^2}
      \stackrel{\eqref{eq6ad89752}}{\le} 1-\frac\sigma4
    =1-\frac{p(s)^2}{2^6{\mu}^2} .
\end{equation}

We shall consider the {\em endpoint
increment} $\Delta \Psi_{\mu}$  of $\Psi_{\mu}$ between the endpoints of $\eta$, equivalently, by \eqref{endpoint1} and \eqref{endpoint2} between the endpoints of $\gamma$:
\[ \Delta \Psi_{\mu} \coloneq \Psi_{\mu}(\eta(L)) -\Psi_{\mu}(\eta(0)) .\]

On the one hand, by \eqref{eq6ad91b1a} we have
\begin{eqnarray}
 \nonumber   \Delta\Psi_{\mu}&\le&
    \int_0^L\abs{\nabla_H\Psi_{\mu}}(\eta(s))\dd s \label{eq6ad91b19}\\
\nonumber    &\stackrel{\eqref{eq6ad91b1a}}{\le}& L -   \int_0^L  \frac{p(s)^2}{2^6{\mu}^2}\dd s\\
    &\stackrel{\eqref{peak_estimate_used}}{\le}& L-\frac{1}{2^6{\mu}^2}\frac{{\mu}^3}{12{\tempo}}
    =L-\frac{{\mu}}{768{\tempo}} .
    \label{eq6ad897d7}
\end{eqnarray}
On the other hand, using the fact that the $ {w}$-coordinate vanishes at both the endpoints of the curve, we bound
\begin{equation}\label{eq6ad897d8}
\begin{aligned}
    \Delta\Psi_{\mu}
    &\overset{\eqref{def_Phi}}=\Delta F_{\mu}-4{\mu}^2\int_{-{\tempo}}^{\tempo}\abs{{\bi}_{\mu}'}^2\dd y\\
    &\overset{\eqref{eq6ad897cd}\&\eqref{eq6ad897c8}}\ge L(\gamma)-\frac92{\mu}-4{\mu}^2\cdot\frac{1225}{\mu}
    >L(\gamma)-5000{\mu} .
\end{aligned}
\end{equation}
Combining~\eqref{eq6ad897d7} and~\eqref{eq6ad897d8},  gives
\[
 L-L(\gamma)
 \ge {\mu}\left(\frac1{768{\tempo}}-5000\right) >0
\]
where the last inequality follows from the choice of ${\tempo}=10^{-8}$.
We reached a contradiction with $L\le L(\gamma)$, which concludes the proof in the first case.

\subsubsection{Case: ${\mu}=0$}
Suppose next that $\mu=0$, so $P\circ\eta\equiv0$, and choose an independent parameter $\varepsilon>0$.
Equations~\eqref{eq6ad897cc} and~\eqref{eq6ad897c68}, applied with this parameter, give
\begin{equation}\label{nabla_uno}
 |\nabla_HF_\varepsilon(\eta(s))|
 =\sqrt{a_\varepsilon(\eta(s))^2+b_\varepsilon(\eta(s))^2}=1.
\end{equation}
Consequently, the endpoint identities~\eqref{endpoint1}--\eqref{endpoint2}, the chain rule, and Cauchy--Schwarz yield
\begin{eqnarray*}
 \Delta F_\varepsilon
 &\stackrel{\eqref{endpoint1},\,\eqref{endpoint2}}{=}&
 \int_0^L\frac{\dd}{\dd s}[F_\varepsilon(\eta(s))]\dd s\\
 &\stackrel{\text{chain rule}}{\le}&
 \int_0^L\abs{\nabla_HF_\varepsilon(\eta(s))\cdot\dot\eta(s)}\dd s\\
 &\stackrel{\text{Cauchy--Schwarz}}{\le}&
 \int_0^L\norm{\nabla_HF_\varepsilon(\eta(s))}\,\norm{\dot\eta(s)}\dd s\\
 &\stackrel{\eqref{nabla_uno},\,\eqref{bound_pbal}}{=}&L.
\end{eqnarray*}
We deduce
\[
 L\ge\Delta F_\varepsilon
 \stackrel{\eqref{eq6ad897cd}}{\ge}L(\gamma)-\frac92\varepsilon.
\]
Letting $\varepsilon\downarrow0$ proves $L\ge L(\gamma)$, and the assumed opposite inequality gives $L=L(\gamma)$.
Together with the case $\mu>0$, this proves minimality of $\gamma$ against every absolutely continuous horizontal competitor.
\end{proof}


\begin{thebibliography}{99}
\small

\bibitem{ABB}
A.~Agrachev, D.~Barilari, U.~Boscain,
\emph{A Comprehensive Introduction to Sub-Riemannian Geometry},
Cambridge Stud.\ Adv.\ Math.\ \textbf{181}, Cambridge University Press, 2020.

\bibitem{BCJPS}
D.~Barilari, Y.~Chitour, F.~Jean, D.~Prandi, M.~Sigalotti,
\emph{On the regularity of abnormal minimizers for rank 2 sub-Riemannian structures},
J.\ Math.\ Pures Appl.\ \textbf{133} (2020), 118--138.

\bibitem{BDFPR}
A.~Belotto da Silva, A.~Figalli, A.~Parusi\'nski, L.~Rifford,
\emph{Strong Sard Conjecture and regularity of singular minimizing geodesics for analytic sub-Riemannian structures in dimension 3},
Invent.\ Math.\ \textbf{229} (2022), 395--448.

\bibitem{CJMRSSS}
Y.~Chitour, F.~Jean, R.~Monti, L.~Rifford, L.~Sacchelli, M.~Sigalotti, A.~Socionovo,
\emph{Not all sub-Riemannian minimizing geodesics are smooth},
preprint, 2025, arXiv:2501.18920.

\bibitem{Fathi}
A.~Fathi,
\emph{Weak KAM Theorem in Lagrangian Dynamics},
preliminary version 10, 2008.

\bibitem{Corners}
E.~Hakavuori, E.~Le Donne,
\emph{Non-minimality of corners in subriemannian geometry},
Invent.\ Math.\ \textbf{206} (2016), 693--704.

\bibitem{Blowups}
E.~Hakavuori, E.~Le Donne,
\emph{Blowups and blowdowns of geodesics in Carnot groups},
J.\ Differential Geom.\ \textbf{123} (2023), no.~2, 267--310.

\bibitem{HL}
R.~Harvey, H.~B.~Lawson, Jr.,
\emph{Calibrated geometries},
Acta Math.\ \textbf{148} (1982), 47--157.

\bibitem{Leonardi_Monti}
G.~P.~Leonardi, R.~Monti,
\emph{End-point equations and regularity of sub-Riemannian geodesics},
Geom.\ Funct.\ Anal.\ \textbf{18} (2008), no.~2, 552--582.

\bibitem{LS}
W.~Liu, H.~J.~Sussmann,
\emph{Shortest paths for sub-Riemannian metrics on rank-two distributions},
Mem.\ Amer.\ Math.\ Soc.\ \textbf{118} (1995), no.~564.

\bibitem{LZreg}
L.~Lokutsievskiy, M.~Zelikin,
\emph{Derivatives of sub-Riemannian geodesics are $L_p$-H\"older continuous},
ESAIM Control Optim.\ Calc.\ Var., 2023, doi:10.1051/cocv/2023055.

\bibitem{Monti}
R.~Monti,
\emph{The regularity problem for sub-Riemannian geodesics},
in: Geometric Control Theory and Sub-Riemannian Geometry, Springer INdAM Ser.\ \textbf{5} (2014), 313--332.

\bibitem{MPV}
R.~Monti, A.~Pigati, D.~Vittone,
\emph{Existence of tangent lines to Carnot--Carath\'eodory geodesics},
Calc.\ Var.\ Partial Differential Equations \textbf{57} (2018), Art.~75.

\bibitem{Palais1957}
R.~S.~Palais,
\emph{A global formulation of the Lie theory of transformation groups},
Mem.\ Amer.\ Math.\ Soc.\ 22, 1957.

\bibitem{RSS}
T.~Rossi, A.~J.~A.~Schiavoni Piazza, A.~Socionovo,
\emph{Interior singularity and branching of geodesics in real-analytic sub-Riemannian manifolds},
preprint, 2026, arXiv:2603.23068.

\bibitem{Soc}
A.~Socionovo,
\emph{Sharp regularity of sub-Riemannian length-minimizing curves},
Calc.\ Var.\ Partial Differential Equations \textbf{65} (2026), Art.~181,
\href{https://doi.org/10.1007/s00526-026-03343-z}{doi:10.1007/s00526-026-03343-z}.

\end{thebibliography}
\end{document}